\documentclass[a4paper]{article}
\pdfoutput=1
\usepackage{subfigure}
\usepackage{graphicx}
\usepackage{lipsum}
\usepackage{amsmath,amssymb,mathrsfs}
\usepackage{hyperref}
\usepackage{url}
\usepackage{mathtools}
\usepackage{changepage}
\usepackage{multirow}
\usepackage{wrapfig}
\usepackage{color}
\usepackage{epsfig,enumerate,amsmath,amsfonts,amssymb,amsthm,mathrsfs,ifpdf}
\usepackage{indentfirst,relsize}
\usepackage{setspace,graphicx}
\usepackage{latexsym}
\usepackage[all]{xy}
\usepackage[usenames,dvipsnames]{pstricks}
\usepackage{pst-grad} 
\usepackage{pst-plot} 
\usepackage[margin = 2.50cm]{geometry}
\usepackage{subcaption}
\usepackage{algorithm}

\usepackage[noend]{algpseudocode}

\newcommand{\remove}[1]{}

\newtheorem{theo}{Theorem}
\newtheorem{lem}[theo]{Lemma}

\newtheorem{coro}[theo]{Corollary}

\newtheorem{cl}[theo]{Claim}

\newtheorem{defi}[theo]{Definition}

\graphicspath{{./images/}}

\title{{Liar's vertex-edge domination in unit disk graph}}

\author{Debojyoti Bhattacharya\footnote{Indian Instittute of Technology Patna, Bihta, 801106, Bihar, India. email: debojyoti\_2021ma11@iitp.ac.in} \and Subhabrata Paul\footnote{Indian Instittute of Technology Patna, Bihta, 801106, Bihar, India. email: subhabrata@iitp.ac.in}}

\date{}

\begin{document}

\maketitle

\begin{abstract}
Let $G=(V, E)$ be a simple undirected graph. A closed neighbourhood of an edge $e=uv$ between two vertices $u$ and $v$ of $G$, denoted by $N_G[e]$, is the set of vertices in the neighbourhood of $u$ and $v$ including $\{u,v\}$. A subset $L$ of $V$ is said to be liar's vertex-edge dominating set if $(i)$ for every edge $e\in E$, $|N_G[e]\cap L|\geq 2$ and $(ii)$ for every pair of distinct edges $e,e'$, $|(N_G[e]\cup N_G[e'])\cap L|\geq 3$. The minimum liar's vertex-edge domination problem is to find the liar's vertex-edge dominating set of minimum cardinality. 
In this article, we show that the liar's vertex-edge domination problem is NP-complete in unit disk graphs, and we design a polynomial time approximation scheme(PTAS) for the minimum liar's vertex-edge domination problem in unit disk graphs.  

\noindent\textbf{keywords:} {Liar's vertex-edge dominating set, NP-completeness, unit disk graphs, approximation scheme}
\end{abstract}

\section{Introduction}

Let $G=(V, E)$ be a graph, where $V$ is the vertex set and $E$ is the edge set, respectively. Also, assume that the set $N_G(v)=\{u\in V|uv\in E\}$ is the open neighbourhood and the set $N_G[v]=N_G(v)\cup \{v\}$ is the closed neighbourhood of a vertex $v\in V$. The closed neighbourhood of an edge $e=xy\in E$ is defined as the set $N_G[e]=N_G[xy]=N_G[x]\cup N_G[y]$.
If $v\in N_G[e]$ for some edge $e$, then we say that $e$ is \emph{vertex-edge dominated} (or \emph{ve-dominated}) by the vertex $v$.   
A \emph{vertex-edge dominating set} (or a \emph{ve-dominating set}) $D$ of $G$ is a subset of the vertex set $V$ such that for every edge $e\in E$, $N_G[e]$ contains at least one vertex of $D$.
The notion of ve-domination was introduced by Peters in his Ph.D. thesis in 1986 \cite{peters}. After that different researchers have studied ve-domination in literature \cite{lewis,boutrig2016vertex,jena,paul2,zylinski2019vertex}.
For some non negative integer $k\geq 2$, an edge $e$ is \emph{$k$-vertex edge dominated}(or \emph{$k$-ve dominated}) by $D_k\subseteq V$ if $D_k$ contains at least $k$ vertices of $N_G[e]$. For some integer $k\geq 2$, A \emph{$k$-ve dominating set} of $G$ is a subset $D_k$ of the vertex set $V$ such that every edge $e\in E$ is $k$-ve dominated by $D_k$. 
$k$-ve domination is also studied in various literatures\cite{krishna,li2023polynomial,naresh, bhattacharya2024k}.
Other variations of ve-domination, namely independent ve-domination \cite{lewis,paul}, total ve-domination \cite{ahangar2021total,Totalve-domchellali}, global ve-domination \cite{chitra2012global,globalvedom} etc are studied in the literature. Liar's vertex-edge domination is a new variation of ve-domination problem that arises in communication networks due to the following circumstances:  

Consider a communication network where each vertex $v\in V$ is a communication device and each edge $e=uv$ is a communication channel between two devices $u$ and $v$. For the flawless functioning of the network, every channel should be constantly monitored by sentinels placed at vertices so that every damaged channel can be identified. A sentinel at some vertex $v\in V$ can detect a damaged channel incident to $N_G[v]$ and then report that channel as damaged. If any sentinel fails to detect a broken or damaged channel, then to correctly identify that damaged channel, each channel is required to be monitored by more than one sentinel. Now, consider the scenario where a sentinel assigned at some vertex $v$ correctly detects the damaged channel incident to a vertex of the neighbourhood of $v$ but misreports an undamaged channel incident to a vertex of the neighbourhood of $v$ as a damaged channel. Also, suppose that at most one sentinel at $v$ in the neighbourhood of a damaged channel $e$ is misreporting an undamaged channel incident to a vertex of the neighbourhood of $v$ as a damaged channel and every other sentinel placed in the vertices of $N_G[e]$ correctly detects and report the damaged or broken channel as $e$. To deal with such problems \emph{liar's vertex-edge dominating set} (or liar's ve-dominating set) was introduced and the following definition was established in \cite{bhattacharya2024liar}.

\begin{defi}
	A subset $L$ of $V$ is said to be a \emph{liar's ve-dominating set} of a graph $G=(V,E)$ if 
	\begin{itemize}
		\item[(i)] for every $e_i\in E$, $|N_G[e_i]\cap L|\geq 2$ and
		\item[(ii)] for every pair of distinct edges $e_i$ and $e_j$, $|(N_G[e_i]\cup N_G[e_j])\cap L|\geq 3$. 
	\end{itemize}
\end{defi}

The minimum liar's ve-domination problem and its corresponding decision version are described as follows:

\noindent\underline{\textsc{Minimum Liar's VE-Domination Problem} (\textsc{MinLVEDP})}

\noindent\emph{Instance}: A graph $G=(V,E)$.

\noindent\emph{Output}: A minimum liar's ve-dominating set $L_{ve}$ of $G$.

\noindent\underline{\textsc{Liar's VE-Domination Decision Problem} (\textsc{DecideLVEDP})}

\noindent\emph{Instance}: A graph $G=(V,E)$ and an integer $k$.

\noindent\emph{Question}: Does there exist a liar's ve-dominating set of size at most $k$?

A graph is a unit disk graph if it is the intersection graph of disks of equal radii in the plane. The most common application of unit disk graphs is in broadcast networks. Let $\{c_1,c_2,\ldots,c_n\}$ be the centres of $n$ circular disks with diameter $1$. The corresponding unit disk graph $G=(V, E)$ is defined as follows: corresponding to every disk centered at $c_i$, there is a vertex $v_i\in V$, and two vertices of $G$ are adjacent if the corresponding disks intersect. Note that, if two disks intersects, then the distance between the center of the two disks is at most $1$. 
This article examines the algorithmic aspects of the liar's ve-domination problem in unit disk graphs.
The contributions are as follows:

\begin{itemize}
	\item[(i)] We show that the \textsc{DecideLVEDP} is NP-complete in unit disk graphs.
	\item[(iii)] We design a PTAS for \textsc{MinLVEDP} in unit disk graphs. 
\end{itemize}

The rest of the paper is organized as follows. In Section $2$, we show that the \textsc{DecideLVEDP} is NP-complete in the unit disk graph. In Section $3$, we propose a PTAS for \textsc{MinLVEDP} in the unit disk graph. Finally, Section $4$ concludes this article with suggestions for future research.

\section{NP-completeness in unit disk graph}\label{sec:NP-completeunitdisc}

In this section, we show that the \textsc{DecideLVEDP} is NP-complete in unit disk graphs. To prove NP-completeness, we describe a polynomial time reduction from an instance of the vertex cover problem for planar graphs with maximum degree $3$, which is known to be NP-complete \cite{garey1979} to an instance of \textsc{DecideLVEDP}. The decision version of the vertex cover problem is defined as follows:

\noindent\underline{\textsc{Vertex Cover Problem in planar graphs} (\textsc{DecideVCPL})}

\noindent\emph{Instance}: A planar graph $G=(V,E)$ with maximum degree $3$ and a positive integer $k$.

\noindent\emph{Question}: Does there exist a set $C\subseteq V$ of size at most $k$ such that for every edge $e=uv$ of $G$ either $u\in C$ or $v\in C$? 

To prove the NP-completeness, we consider a special type of embedding of the planar graph $G$. The following lemma describes the embedding of $G$.
\begin{lem}[\cite{embeddingofplanerg}]\label{lem:planarembedding}
	A planar graph $G=(V,E)$ with maximum degree $4$ can be embedded in the plane using $O(|V|^2)$ area in such a way that its vertices are at integer coordinates and its edges are drawn so that they are made up of line segments of the form $x=i$ or $y=j$, for integers $i$ and $j$.
\end{lem}

Such an embedding is called the orthogonal embedding of a graph. A linear time algorithm to find an orthogonal drawing of a given graph with at most $2$ bends along each edge is given in \cite{BIEDL1998159}.

\begin{coro}\label{coro:planarembedding}
	A planar graph $G=(V, E)$ with maximum degree $3$ and $|V|\geq 3$ can be embedded in the plane with its vertices at $(10i,10j)$ and its edges are drawn as a sequence of consecutive axis-parallel line segments on the lines $x=10i$ or $y=10j$, for some $i$ and $j$.
\end{coro}

\begin{theo}
	The \textsc{DecideLVEDP} is NP-complete in unit disk graphs.
\end{theo}

\begin{proof}
	Observe that, \textsc{DecideLVEDP} is in $NP$. Next, we show a polynomial time reduction of an instance $G=(V, E)$ of \textsc{DecideVCPL} to an instance $G'=(V', E')$ of \textsc{DecideLVEDP} which is a unit disk graph.

		First, consider an embedding of the planar graph $G$. The embedding can be constructed using one of the algorithms in \cite{Efficientplanarity,hamiltonianpathingreedgraph}. 
		The construction of the graph $G'$ is as follows:
		\begin{figure}
			\centering
			\includegraphics[scale=0.5]{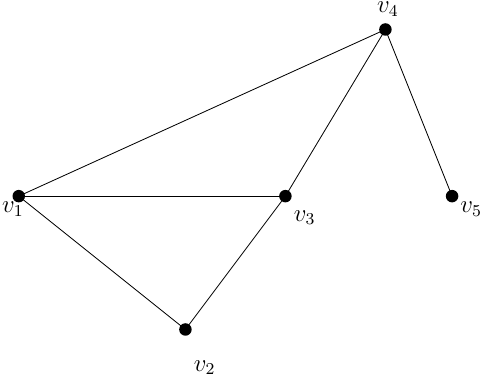}
			\caption{A planar graph $G$ with maximum degree $3$}
			\label{fig:planargraph}
		\end{figure}
		
		\begin{figure}
			\centering
			\includegraphics[scale=0.5]{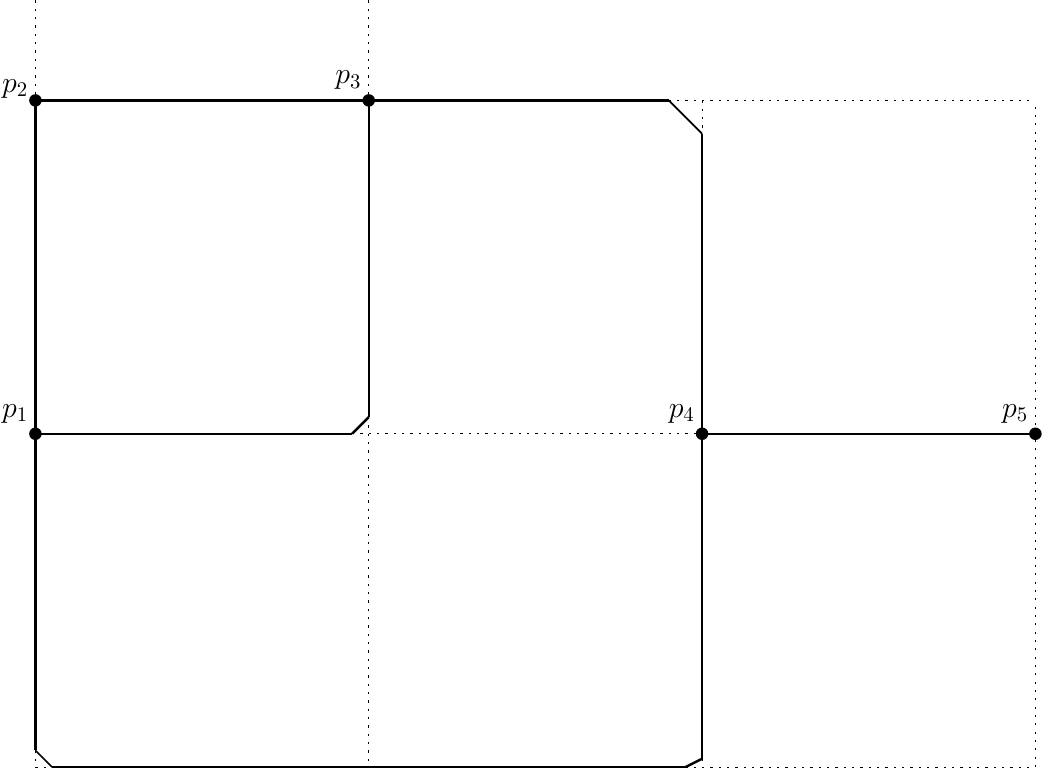}
			\caption{Embedding of the graph $G$ on a grid.}
			\label{fig:planarembedding}
		\end{figure}

		\begin{itemize}
			\item[(i)] Using one of the algorithms in \cite{Algorithmforembedding1,Algorithmforembedding2} and Lemma \ref{lem:planarembedding} and Corollary \ref{coro:planarembedding}, $G$ can be embedded in plane such that each edge of $G$ is a sequence of connected line segments of length $10$ units. Let $l$ be the total number of line segments in the embedding and the points $P=\{p_1,p_2,\ldots,p_n\}$ be the node points corresponding to the vertices $V=\{v_1,v_2,\ldots,v_n\}$.
			\item[(ii)]Note that in the embedding, each edge is a sequence of line segments. To construct a unit disk graph, we add some points in the line segments such that the distance between any consecutive points is at most $1$. In other words, the points added in the line segments joining $p_i$ and $p_j$ forms a path between $p_i$ and $p_j$.Now, we describe an efficient way to add the points in the line segments. For every edge $p_ip_j$ of length $10$, add a point at a distance $0.8$ from $p_i$ and add another ten points at a distance $1.2, 2.2, 3.2, 4.2, 4.8,5.3,6.3,7.3,8.3,$ and $9.2$ from $p_i$, respectively. For every other edge of $G$, the embedding has more than one line segment between $p_i$ and $p_j$. We add some points in the line segments in two steps. First, we add eleven points to any one of the segments incident to a node point at a distance $1,1.5,2,3,4,4.9,5.9,$ $6.9,7.9,8.9,$ and $9.9 $ from the node point. For the other segment incident to another node point add ten points at a distance $1,2,3,4,4.9,5.9,$ $6.9,7.9,8.9,$ and $9.9$ from that node point. Finally, add ten points to the remaining line segments at a distance of $0.5,1.5,2.5,3.5,4.5,5.5,6.5,7.5,$ $8.5$, and $9.5$ from any one of the endpoints of the line segment. Let $A$ be the set of points added in this step.
			\item[(iii)] Add a line segment of length $4$ to every $p_i\in P$. This is possible because the grid is a $10\times 10$ grid and the maximum degree of $G$ is $3$. In the line segment add four points at a distance $1,2,3,$ and $4$ from $p_i\in P$. Let $S$ be the set of points added in this step.
		\end{itemize}
		
		\begin{figure}
			\centering
			\includegraphics[scale=0.5]{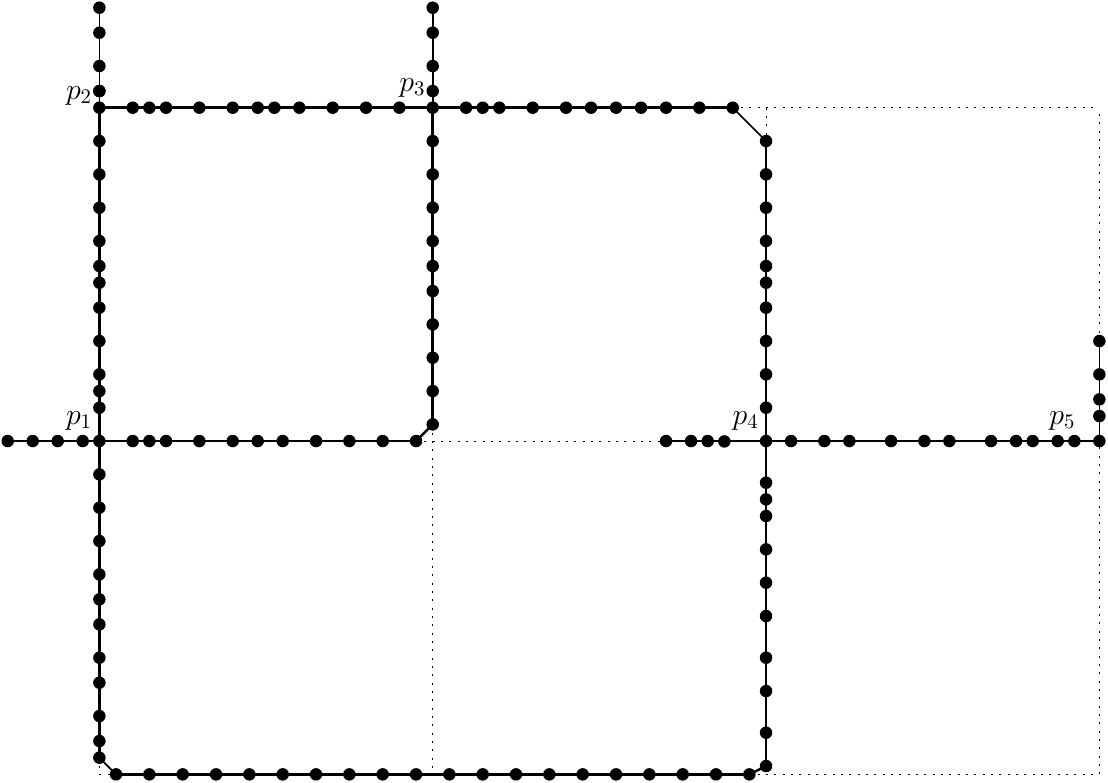}
			\caption{Construction of unit disk graph $G'$ from embedding of $G$}.
			\label{fig:construction}
		\end{figure}
		The total number of points in $A$ is $10l+|E|$. Let $A=\{q_1,q_2,\ldots,q_{10l+|E|}\}$ and $S=\{u_i,x_i,y_i,z_i:1\leq i\leq n\}$. Consider the graph $G'=(V',E')$, where $V'=P\cup A \cup S$. Note that, two vertices of $G'$ are adjacent if and only if the Euclidean distance between the points is at most $1$. Thus, the graph $G'$ is a unit disk graph. Also, $|V'|=|P|+10l+|E|+3n=4n+10l+|E|$. Since $G$ is planar, it is easy to see that $G'$ can be constructed in polynomial time. An example of a planar graph $G$ with maximum degree $3$ is given in Figure \ref{fig:planargraph}. The embedding of $G$ is given in Figure \ref{fig:planarembedding}. The construction of the corresponding unit disk graph $G'$ from the embedding of $G$ is given in Figure \ref{fig:construction}.
	
	Next, we prove the following lemma.
	
	\begin{lem}
		$G$ has a vertex cover of size at most $k$ if and only if $G'$ has a liar's ve-dominating set of size at most $k+3n+6l$.
	\end{lem}
	
	\begin{proof}
		Let $C$ be a cover of $G$ of size at most $k$. Let $L'=C\bigcup\limits_{i=1}^{n} \{u_i,x_i,y_i\}$. We construct a set $A'$ containing the points of $A$ by choosing some vertices from each segment. For any edge $v_iv_j$ of $G$, at least one of the vertices corresponding to $\{v_i,v_j\}$ is in $C$. Without loss of generality, let us assume that $v_i\in C$. Let $p_i$ and $p_j$ be the vertices of $G'$ corresponding to the vertices of $v_i$ and $v_j$, respectively. 
		
		\begin{figure*}[ht!]
		
				\begin{subfigure}[]
					{
					\includegraphics[scale=0.45]{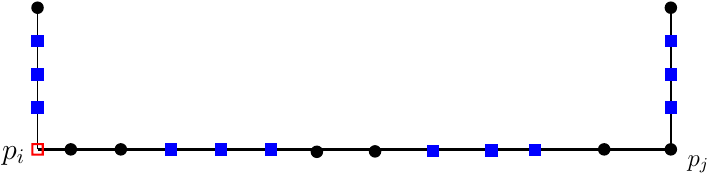}
					\label{fig:singlesegment}}
				\end{subfigure}%
				~
				\begin{subfigure}[]{
					\includegraphics[scale=0.42]{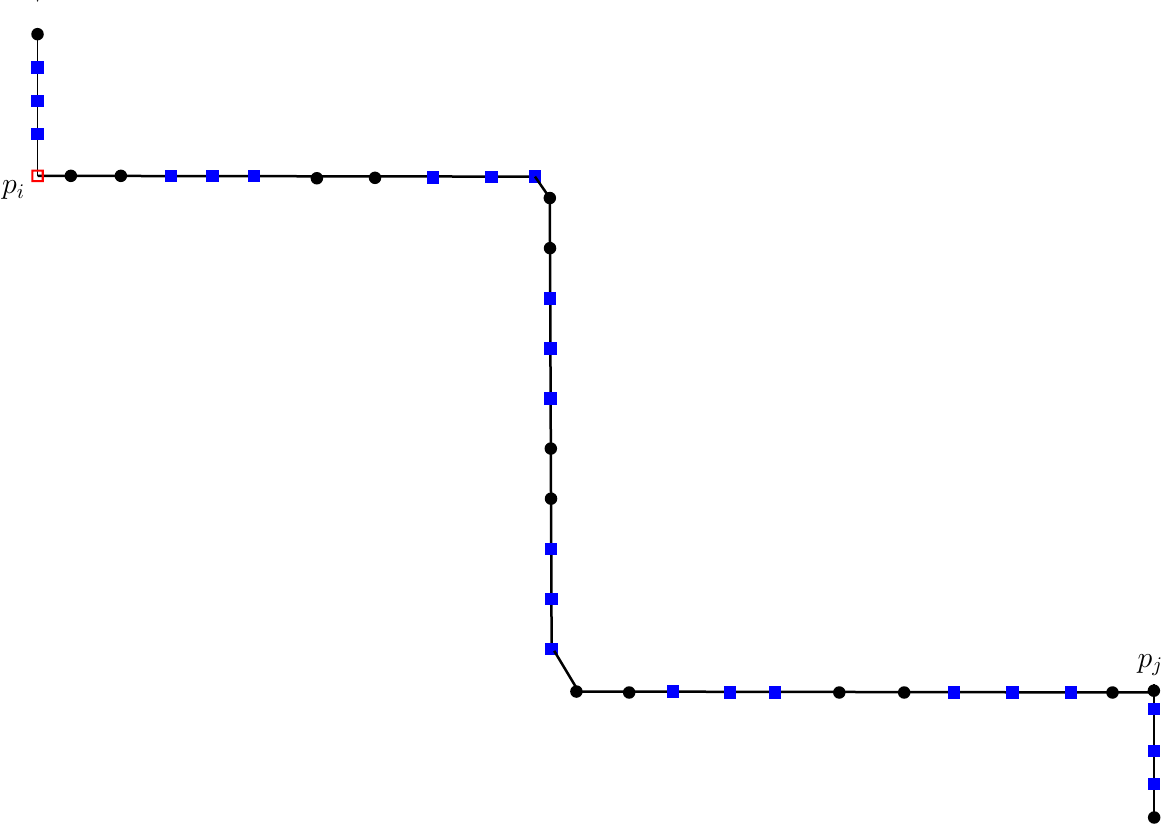}
					\label{fig:multiplesegment}}
				\end{subfigure}
				\caption{In this figure, we demonstrate the vertices to be selected in $L$. Suppose $v_i\in C$. Corresponding $p_i$ is marked as a box vertex. The square vertices are the vertices from the segments incident to $p_i$ and $p_j$ which are selected in $L$.}
				\label{fig:lvedset}
			
		\end{figure*}

		For every edge which is represented by a single line segment of length $10$, suppose that $\{p_i,q_{i+1},q_{i+2},q_{i+3}\\,q_{i+4},q_{i+5},q_{i+6},q_{i+7},q_{i+8},q_{i+9},q_{i+10},q_{i+11},p_j\}$ is the path in $G'$ corresponding to the segment $p_ip_j$. We include the vertices $\{q_{i+3},q_{i+4},q_{i+5},q_{i+8},q_{i+9},q_{i+10}\}$ in $A'$(see Fig \ref{fig:singlesegment}). For every edge represented by a sequence of line segments, there are two cases. Note that, from the construction of $G'$, for every such edge, one of the line segments incident to a node point has eleven added points and the other line segment incident to another node point has ten added points. First, without loss of generality, we assume that the segment incident to $p_i$ has eleven added points. Add the third, fourth, fifth, eighth, ninth, and tenth points from $p_i$, and in the following segments select the second, third, fourth, seventh, eighth, and ninth added points in $A'$. In the second case (the segment incident to $p_i$ has ten vertices), we select the third, fourth, fifth, eighth, ninth, and tenth points from $p_i$ in $A'$(see Figure \ref{fig:multiplesegment}). In the following segments, we select three consecutive vertices after skipping two consecutive vertices starting with skipping the first two vertices of the segment adjacent to the segment incident to $p_i$. Observe that, proceeding in this way, every segment has six vertices in $A'$. Therefore, if an edge is represented by $l'$ line segments, then the total number of points in $A'$ corresponding to that edge is $6l'$. We follow this procedure for every edge in $G$. Thus, we have $|A'|=6l$. Now, let $L=L'\cup A'$. It is easy to verify that $L$ is a liar's ve-dominating set of $G'$ and $|L|\leq k+3n+6l$.
		
		Conversely, suppose that $L$ is a liar's ve-dominating set of $G'$ of size at most $k+3n+6l$. For every edge $y_iz_i$ and $x_iy_i$, we have $|(N_{G'}[y_iz_i]\cup N_{G'}[x_iy_i])\cap L|\geq 3$. If all four vertices of $\{u_i,x_i,y_i,z_i\}$ are in $L$, then we can obtain a liar's ve-dominating set of $G'$ of size at most $k+3n+6l$ by removing $z_i$ from $L$. If one of the vertices in $\{u_i,x_i,y_i\}$ is not in $L$, then obviously $z_i\in L$ and we can always include the missing vertex in $L$ by replacing $z_i$. Therefore, without loss of generality, we assume that $\{u_i,x_i,y_i\}$ is in $L$ for every $p_i$. Now, let $L'=L\setminus(\bigcup\limits_{i=1}^{n}\{u_i,x_i,y_i\})$. 
		\begin{cl}\label{claim:atleast6vertex}
			Every segment has at least six vertices in $L$.
		\end{cl}
		\begin{proof}
			Consider a segment of length $10$, say $\{p_i,q_{i+1},q_{i+2},q_{i+3},q_{i+4},q_{i+5},q_{i+6},\\q_{i+7}$  $,q_{i+8},q_{i+9},q_{i+10},q_{i+11},p_j\}$, which is incident to the node points $p_i$ and $p_j$. To satisfy condition $(2)$ of liar's ve-domination, there must be at least three vertices of $N_{G'}[q_{i+1}q_{i+2}]\cup N_{G'}[q_{i+2}q_{i+3}]$ in $L$. Now, $N_{G'}[q_{i+1}q_{i+2}]\cup N_{G'}[q_{i+2}q_{i+3}]=\{p_i,q_{i+1},q_{i+2},q_{i+3},q_{i+4}\}$. Hence, $L$ contains at least three vertices of the set $\{p_i,q_{i+1},q_{i+2},q_{i+3},q_{i+4}\}$. Also, for the pair of edges $q_{i+6}q_{i+7}$ and $q_{i+7}q_{i+8}$, there must be at least three vertices of $N_{G'}[q_{i+6}q_{i+7}]\cup N_{G'}[q_{i+7}q_{i+8}]$ in $L$. Thus, $L$ contains at least three vertices of $\{q_{i+5},q_{i+6},q_{i+7},q_{i+8},q_{i+9}\}$.		
			This implies that at least six vertices of the segment must be in $L$.
			
			Consider an edge $v_iv_j$ of $G$ which is represented by more than one line segment in $G'$. Let $p_i$ and $p_j$ be the corresponding points of $v_i$ and $v_j$ in $G'$.  Let us assume that there are $l'$ line segments joining $p_i$ and $p_j$. Now, each line segment except one in $G'$ has ten vertices, and only one line segment which is either incident to $p_i$ or $p_j$ has eleven vertices.  Now, we show that for every edge represented by $l'$ number of line segments, there are at least $6l'$ vertices in $L$.
			\begin{figure}
				\centering
				\includegraphics[scale=0.6]{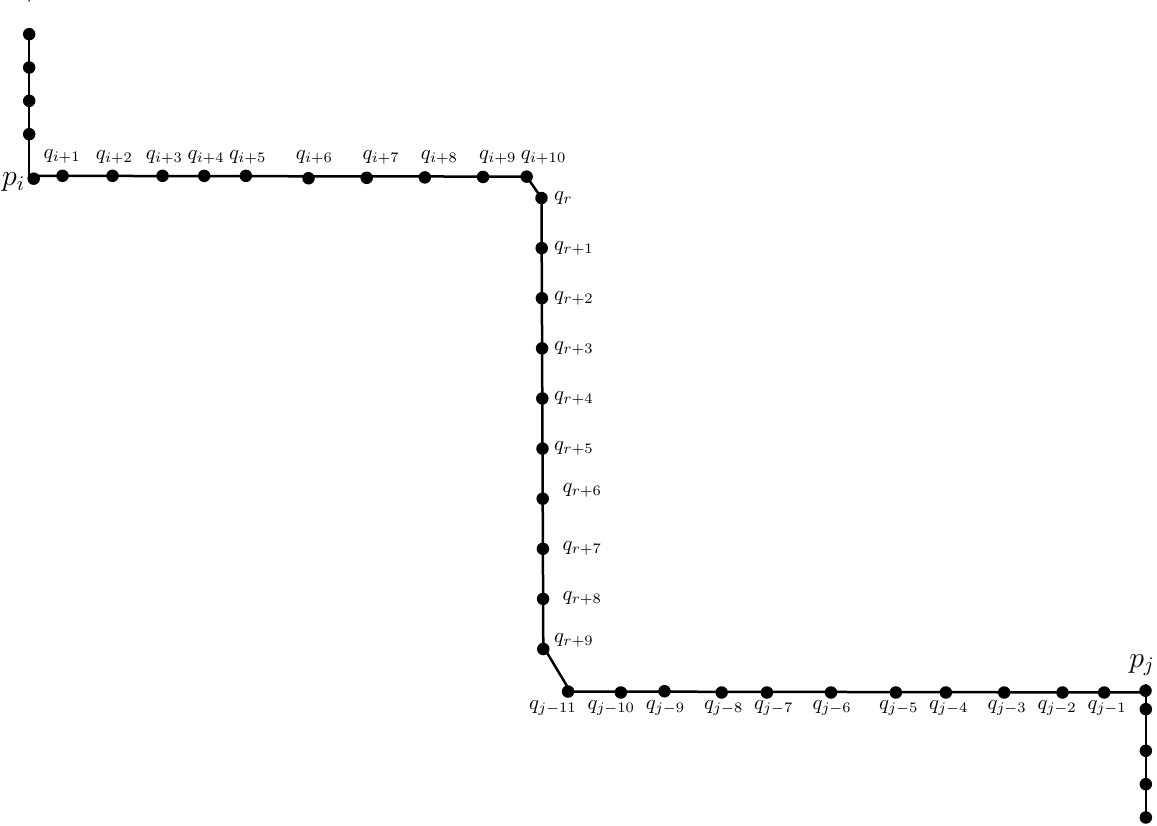}
				\caption{A path between $p_i$ and $p_j$ where every segment has at least $6$ vertices in $L$}
			\end{figure}
			Let us assume that there are ten vertices $\{q_{i+1},q_{i+2},q_{i+3},q_{i+4},q_{i+5},q_{i+6},q_{i+7},$ $q_{i+8},q_{i+9},q_{i+10}\}$ in the line segment incident to $p_i$. For the edges $q_{i+2}q_{i+3}$ and $q_{i+3}q_{i+4}$, there must be at least three vertices of $\{q_{i+1},q_{i+2},q_{i+3},q_{i+4},q_{i+5}\}$ is in $L$. Also, for the pair of edges $q_{i+6}q_{i+7}$ and $q_{i+7}q_{i+8}$,  there must be at least three vertices of $\{q_{i+5},q_{i+6},q_{i+7},q_{i+8},q_{i+9}\}$ is in $L$. Hence, $L$ contains six vertices corresponding to the line segment incident to $p_i$.
			
			Now, consider a line segment that is not incident to any of the vertices in $p_i$ and $p_j$. According to the construction of $G'$, there are ten vertices $\{q_r,q_{r+1},q_{r+2},$ $q_{r+3}, q_{r+4},q_{r+5},q_{r+6},q_{r+7},q_{r+8}, q_{r+9}\}$ in the segment. Observe that, to satisfy the condition $(2)$ of the liar's ve-domination, for the edges $q_{r+1}q_{r+2}$ and $q_{r+2}q_{r+3}$ there must be at least three vertices of $\{q_r,q_{r+1},q_{r+2},$ $q_{r+3},q_{r+4}\}$ is in $L$. Also, for the pair of edges $q_{i+6}q_{i+7}$ and $q_{i+7}q_{i+8}$,  there must be at least three vertices of $\{q_{i+5},q_{i+6},q_{i+7},q_{i+8},q_{i+9}\}$ is in $L$. Hence, $L$ contains six vertices corresponding to such the line segment.
			
			From the construction, segment incident to $p_j$ must have eleven vertices $\{q_{j-1},q_{j-2},$ $q_{j-3},q_{j-4},q_{j-5},\\q_{j-6},q_{j-7},q_{j-8},q_{j-9},q_{j-10},q_{j-11}\}$, where $q_{j-1}p_j$ is an edge of $G'$. Now, for the edges $q_{j-2}q_{j-3}$ and $q_{j-3}q_{j-4}$, there must be at least three vertices of $\{q_{j-1},q_{j-2},q_{j-3},q_{j-4},q_{j-5}\}$ in $L$. Also, for the pair of edges $q_{j-7}q_{j-8}$ and $q_{j-8}q_{j-9}$,  there must be at least three vertices of $\{q_{j-6},q_{j-7},q_{j-8},q_{j-9},q_{j-10}\}$ is in $L$. Hence,  $L$ contains six vertices corresponding to such the line segment incident to $p_j$.	
		\end{proof}

		Thus, for every line segment, there are at least six vertices in $L$. Hence, the total number of vertices in $L$ corresponding to the edge is $6l'$. Similarly, it can be shown that if the segment incident to $p_i$ has eleven vertices, then there are at least $6l'$ vertices in $L$ corresponding to the edge $p_ip_j$. 
		\begin{cl}\label{cl:6l'+1}
			For every edge $v_iv_j$ in $G$, if the corresponding vertices $p_i$ and $p_j$ in $G'$ are not in $L$. Then there are $6l'+1$ points in $L$ corresponding to the line segments joining $p_i$ and $p_j$, where $l'$ is the number of line segments joining $p_i$ and $p_j$. 
		\end{cl}
		\begin{proof}
			Now, suppose there is an edge $v_iv_j$ in $G$ for which none of the corresponding vertices($p_i$ and $p_j$) of $v_i$ and $v_j$ is in $L$. We claim that in such cases the line segments joining $p_i$ and $p_j$ must have at least $6l'+1$ vertices in $L$, where $l'$ is the number of line segments joining $p_i$ and $p_j$. We have already proved that each segment contributes at least six vertices in $L$.
			Now, let us assume that the segment incident to $p_i$ and $p_j$. Thus, it has eleven vertices, namely $\{q_{i+1},q_{i+2},q_{i+3},q_{i+4},q_{i+5},q_{i+7},q_{i+8},q_{i+9},q_{i+10},q_{i+11}\}$, where $p_i$ is adjacent to $q_{i+1}$ and $p_{i+11}$ is adjacent to $p_j$. To satisfy condition $(2)$ of liar's ve-domination, for the pair of edges  $q_{i+1}q_{i+2}$ and $q_{i+2}q_{i+3}$, there must be at least three vertices in $L$. Hence, $L$ contains at least three vertices of the set $\{q_{i+1},q_{i+2},q_{i+3},q_{i+4}\}$. Similarly, for the pair of edges $q_{i+6}q_{i+7}$ and $q_{i+7}q_{i+8}$, $L$ contains at least three vertices of $\{q_{i+5},q_{i+6},q_{i+7},q_{i+8},q_{i+9}\}$. Now, the neighbourhood of the edge $q_{i+10}q_{i+11}$ contains at most one vertex of the set $\{q_{i+5},q_{i+6},$ $q_{i+7},q_{i+8},q_{i+9}\}$. Thus, to $2$-ve dominate the edge $q_{i+10}q_{i+11}$, at least one  vertex of $\{q_{i+10},q_{i+11},p_j\}$ is in $L$. Since $p_j\notin L$, there must be a vertex of $\{q_{i+10},q_{i+11}\}$ in $L$. This implies that at least seven vertices corresponding to the segment incident to $p_i$ and $p_j$ are in $L$.

			Now, assume that there are $l'$ line segments joining $p_i$ and $p_j$, where $l'>1$ and the segment incident to $p_i$ has eleven vertices. Thus, there are ten vertices in the segment incident to $p_j$ and there are exactly $6l'$ vertices in $L$ corresponding to the line segments joining $p_i$ and $p_j$.
			Let $\{p_i, q_{i+1},q_{i+2},q_{i+3},q_{i+4},q_{i+5},$ $q_{i+7},q_{i+8},q_{i+9},q_{i+10},q_{i+11}\}$ be the set of vertices in the line segment incident to $p_i$. Since $p_i$ is not in $L$, to satisfy the condition $(2)$ of liar's ve-domination for the pair of edges $(q_{i+1}q_{i+2},q_{i+2}q_{i+3})$, there must be three vertices in $L$ from $\{q_{i+1},q_{i+2},q_{i+3},q_{i+4}\}$. If $\{q_{i+1},q_{i+2},q_{i+3}\}$ is in $L$, then to satisfy the condition $(2)$ of liar's ve-domination for the pair of edges $q_{i+5}q_{i+6}, q_{i+6}q_{i+7}$, $L$ must contain three vertices from $\{q_{i+4}, q_{i+5},q_{i+6},q_{i+7},q_{i+8}\}$. Therefore, for the edge $q_{i+9}, q_{i+10}$, $L$ must contain a vertex in $\{q_{i+9},q_{i+10},q_{i+11}\}$. Hence, there are seven vertices in $L$ from the segment incident to $p_i$. Similarly, one can show that if $\{q_{i+1},q_{i+3},q_{i+4}\}$ is in $L$ or $\{q_{i+1},q_{i+3},q_{i+4}\}$ is in $L$, then there are seven vertices in $L$ from the line segment incident to $p_i$.  Therefore, we have $\{q_{i+2},q_{i+3},q_{i+4}\}$ is in $L$ and for the pair of edges $q_{i+6}q_{i+7}$ and $q_{i+7}q_{i+8}$, there are three vertices of the set $\{q_{i+5},q_{i+6},q_{i+7},q_{i+8},q_{i+9}\}$ is in $L$. Since $q_{i+4}$ is in $L$, to $2$-ve dominate the edge $q_{i+5}q_{i+6}$, there must be a vertex of $\{q_{i+5},q_{i+6},q_{i+7}\}$ is in $L$. If $q_{i+5}$ is in $L$, then to satisfy condition $(2)$ of liar's ve-domination at least three vertices of $\{q_{i+6},q_{i+7},q_{i+8},q_{i+9}.q_{i+10}\}$ is in $L$. Thus, there are seven vertices in $L$ from the line segment incident to $p_i$. Similarly, it can be shown that if $q_{i+6}$ is in $L$, then $L$ contains seven vertices from the segment incident to $p_i$.
			Therefore, $L$ contains $q_{i+7}$. Hence, $\{q_{i+7},q_{i+8},q_{i+9}\}$ is in $L$. Note that, for the edge $q_{i+10}q_{i+11}$, there must be $2$ vertices in $L$ from the set $\{q_{i+9},q_{i+10},q_{i+11},q_{i+12}\}$, where $q_{i+12}$ is a vertex of the segment closest to the segment incident to $p_i$. Since $q_{i+9}$ is in $L$, if $q_{i+10}$ or $q_{i+11}$ is in $L$, then there are seven vertices in the segment incident to $p_i$. Otherwise, $q_{i+12}$ is in $L$. Since the segment containing $q_{i+12}$ must have $6$ vertices in $L$, it can be shown that $L$ contains the vertices $q_{i+13},q_{i+14}, q_{i+17},q_{i+18}, q_{i+19}$ from that segment. Hence, the edge $q_{i+20}q_{i+21}$ must be ve-dominated by a vertex of the segment incident to $q_{i+21}$. Thus except for the segment incident to $p_i$, $L$ contains a vertex from the segment that is nearest(in terms of Eucledian distance) to its' previous segment. Proceeding in this way, observe that for the segment incident to $p_j$, $L$ contains the vertex from the segment that is nearest to its' previous segment. Suppose that $\{q_{j-10}, q_{j-9}, q_{j-8}, q_{j-7}, q_{j-6}, q_{j-5}, q_{j-4}, q_{j-3}, q_{j-2}, q_{j-1}, p_j\}$ be the vertices in the line segment incident to $p_j$ and $q_{j-10}$ is the nearest vertex from the segment incident to $p_j$ to its' previous segment. Thus, $q_{j-10}$ is in $L$. Therefore, $L$ contains another five vertices $\{q_{j-9},q_{j-8},q_{j-5},q_{j-4},q_{j-3}\}$ from the segment incident to $p_j$. Now, to $2$-ve dominate the edge $q_{j-2}q_{j-1}$, $L$ must contain a vertex from $\{q_{j-2},q_{j-1},p_j\}$. Hence, the segment incident to $p_j$ must have seven vertices in $L$.

			Therefore, if both $p_i$ and $p_j$ are not in $L$, then there are $6l'+1$ number of vertices of the line segments connecting $p_i$ and $p_j$ in $L$, where $l'$ is the number of segments connecting $p_i$ and $p_j$.	
		\end{proof}
		For the edges corresponding to Claim \ref{cl:6l'+1}, replace the extra vertex with one of the vertices in $p_j$ or $p_j$, and we get another liar's ve-dominating set of size at most $k+3n+6l$. Let $A''$ be the number of vertices in added points which are in $L$. Consider the set $C=L'\setminus A''$. Clearly, $|C|\leq k$. Now, for every edge $p_ip_j$, either $p_i\in C$ or $p_j\in C$. Considering the corresponding vertices of such $p_i$ or $p_j$ in $G$ yields a vertex cover of $G$ of size at most $k$.
	\end{proof}
	
	From the above Lemma, we can conclude that \textsc{DecideLVEDP} is NP-complete for Unit disk graphs.
\end{proof}

\section{PTAS for unit disk graph}
In this section, we present a PTAS for the \textsc{MinLVEDP} problem in unit disk graphs. Given any $\epsilon>0$, a PTAS for the \textsc{MinLVEDP} problem in unit disk graphs is a $(1+\epsilon)$-approximation algorithm. We define the distance $\delta(u,v)$ between any two vertices $u,v\in V$ to be the number of edges in the shortest path between the vertices $u$ and $v$. For any two subsets $V_i$ and $V_j$ of $V$, the distance between $V_i$ and $V_j$ is defined as $\delta(V_i,V_j)=\min\limits_{u\in V_i,v\in V_j}\{\delta(u,v)\}$. Given any subset $V_i$ of $V$, let $L(V_i)$ and $L^*(V_i)$ be the liar's ve-dominating set and optimal liar's ve-dominating set of the set of edges of the subgraph induced by $V_i$ in $G$, respectively. The closed neighbourhood of a subset $V_i$ of $V$ is defined by $N_G[V_i]=\bigcup\limits_{v\in V_i}N_G[v]$. The $r$-th neighbourhood of a vertex $v\in V$ is $N_G^r[v]=\{u:\delta(u,v)\leq r\}$. For our algorithm, we need a collection of disjoint subsets of vertices called $m$-separated collection. The $m$-separated collection is defined as follows:

\begin{defi}[$m$-separated collection]
	Let $\mathcal{S}=\{\mathcal{S}_1,\mathcal{S}_2,\ldots,\mathcal{S}_k\}$ be a collection disjoint subsets of $V$, where $\mathcal{S}_i\subseteq V$ for every $i=1,2,\ldots,k$. $\mathcal{S}$ is a $m$-separated collection of subsets of $V$ if $\delta(\mathcal{S}_i,\mathcal{S}_j)>m$ for every $1\leq i,j\leq k$ and $i\neq j$.   
\end{defi}

Nieberg and Hurink considered a $2$-separated collection to design a PTAS for minimum dominating set problem in unit disk graph \cite{Unitdiscdom}. Jena and Das used the $4$-separated collection to propose a PTAS for minimum ve-domination problem in unit disk graphs \cite{jena}. For the liar's ve-domination problem, we consider a $2$-separated collection. Now, we discuss two lemmas regarding the bound on the optimal liar's ve-dominating set induced by the sets in the $m$-separated collection. Note that, these lemmas are true for any undirected graph.

\begin{lem}\label{lem:unitdiscfirst}
	If $\mathcal{S}=\{\mathcal{S}_1,\mathcal{S}_2,\ldots,\mathcal{S}_k\}$ is a $m$-separated collection in $G=(V,E)$, then $\sum_{i=1}^{k}|L^*(\mathcal{S}_i)|\leq |L^*(V)|$ for any $m\geq 2$.
\end{lem}

\begin{proof}
	For each $\mathcal{S}_i\in \mathcal{S}$, let $\mathcal{T}_i=\{u:v\in \mathcal{S}_i~ and ~\delta(u,v)\leq 1\}$. Clearly, $\mathcal{S}_i\subseteq \mathcal{T}_i$ for every $i=1,2,\ldots,k$. Since $\mathcal{S}$ is $m$-separated and $m\geq 2$, we have $\mathcal{T}_i\cap \mathcal{T}_j=\emptyset$ for every $1\leq i,j\leq k$ and $i\neq j$. This implies that $(\mathcal{T}_i\cap L^*(V))\cap (\mathcal{T}_j\cap L^*(V))=\emptyset$, for every $1\leq i,j\leq k$ and $i\neq j$. Thus, $\sum\limits_{i=1}^{k}|\mathcal{T}_i\cap L^*(V)|\leq |L^*(V)|$. Also, $\mathcal{T}_i\cap L^*(V)$ is a liar's ve-dominating set of $\mathcal{S}_i$ for every $i=1,2,\ldots,k$. Therefore, we have $|L^*(\mathcal{S}_i)|\leq |\mathcal{T}_i\cap L^*(V)|$ as $L^*(\mathcal{S}_i)$ is a minimum liar's ve-dominating set of $\mathcal{S}_i$. Hence, we have $\sum\limits_{i=1}^{k}|L^*(\mathcal{S}_i)|\leq \sum\limits_{i=1}^{k}|\mathcal{T}_i\cap L^*(V)|$. Therefore, we have $\sum\limits_{i=1}^{k}|L^*(\mathcal{S}_i)|\leq |L^*(V)|$.
\end{proof}

\begin{lem}\label{lem:unitdiscsecond}
	Let $\mathcal{S}=\{\mathcal{S}_1,\mathcal{S}_2,\ldots,\mathcal{S}_k\}$ be a $m$-separated collection in $G=(V,E)$, $m\geq 2$ and $\mathcal{T}=\{\mathcal{T}_1,\mathcal{T}_2,\ldots,\mathcal{T}_k\}$ be subsets of $V$ with $\mathcal{S}_i\subseteq \mathcal{T}_i$ for every $i=1,2,\ldots,k$. If there exists a $\rho\geq 1$ such that $$|L^*(\mathcal{T}_i)|\leq \rho |L^*(\mathcal{S}_i)|$$ holds for every $i=1,2,\ldots,k$ and if $\bigcup\limits_{i=1}^{k}L^*(\mathcal{T}_i)$ is a liar's ve-dominating set in $G$, then $\sum\limits_{i=1}^{k}|L^*(\mathcal{T}_i)|$ is a $\rho$-approximation of a minimum liar's ve-dominating set in $G$. 
\end{lem}

\begin{proof}
	From Lemma \ref{lem:unitdiscfirst}, we have $\sum\limits_{i=1}^{k}|L^*(\mathcal{S}_i)|\leq |L^*(V)|$. The following inequality completes the proof: $$\sum\limits_{i=1}^{k}|L^*(\mathcal{T}_i)|\leq \rho \sum\limits_{i=1}^{k}|L^*(\mathcal{S}_i)|\leq \rho |L^*(V)|$$
\end{proof}

Now, we give a brief description of how we proceed to the PTAS. First we find a $2$-separated collection $\mathcal{S}=\{\mathcal{S}_1,\mathcal{S}_2,\ldots,\mathcal{S}_k\}$ and it's corresponding set $\mathcal{T}=\{\mathcal{T}_1,\mathcal{T}_2,\ldots,\mathcal{T}_k\}$. After that, we show that there is a $\rho\geq 1$ such that $|L^*(\mathcal{T}_i)|\leq \rho |L^*(\mathcal{S}_i)|$ for every $i=1,,\ldots,k$. Finally, we prove that $\bigcup\limits_{i=1}^{k}|L^*(\mathcal{T}_i)|$ is a liar's ve-dominating set in $G$ and hence the PTAS follows from Lemma \ref{lem:unitdiscsecond}. 

Next, we describe the process to construct the desired a $2$-separated collection $\mathcal{S}=\{\mathcal{S}_1,\mathcal{S}_2,\ldots,\mathcal{S}_k\}$ and its' corresponding set $\mathcal{T}=\{\mathcal{T}_1,\mathcal{T}_2,\ldots,\mathcal{T}_k\}$ such that $\mathcal{S}_i\subseteq\mathcal{T}_i$, for every $i=1,2,\ldots,k$. The process is iterative. First, we set $V_0=V$ and $i=0$. We start with a vertex $v\in V_i$ and compute the liar's ve-dominating set of the subgraphs induced by $N_G^r[v]$ and $N_G^{r+6}[v]$. The algorithm to compute liar's ve-dominating set of $N_G^r[v]$ is given in Algorithm \ref{Algo:liarindepent}. The algorithm is simple. We compute maximal independent sets $I_i$ for the vertex set $V_i$ and update $V_{i+1}$ by $V_i\setminus I_i$ for $i=1,2,3$ and we consider the set $L=I_1\cup I_2\cup I_3$. In the next lemma, we show that $L$ is a lair's ve-dominating set of $G$. 
\begin{algorithm}
	\textbf{Input:} A graph $G=(V,E)$.\\
	\textbf{Outpu:} Liar's ve-dominating set $L$.
	\caption{Liar's\_ve-domination}
	\label{Algo:liarindepent}
	\begin{algorithmic}[1]
		\State $L=\emptyset$ and $V_1=V$;
		\State $i=0$ and $I_0=\emptyset$;
		\For{$i=1$ to $3$}
		\State Compute a maximal independent set $I_i$ in $V_i$;
		\State $V_{i+1}=V_i\setminus I_i$;
		\State $L=L\cup I_i$; 
		\EndFor
		\State \Return $L$;
	\end{algorithmic}
\end{algorithm}
\begin{lem}\label{lem:liarindependent}
	Algorithm \ref{Algo:liarindepent} outputs a liar's ve-dominating set $L$ of $G=(V,E)$.
\end{lem}
\begin{proof}
	Observe that every edge $e\in E$ is ve-dominated by $I_1$ otherwise, if there is an edge $e'=u'v'$ such that $(N_G[u']\cup N_G[v'])\cap I_1=\emptyset$. This implies that $I_1$ is not the maximal independent set. Thus, every edge in $E$ is ve-dominated by $I_1$. Note that, both endpoints of $e$ must have a neighbour in $I_1$. By the same argument as above, every edge in $G'=(V',E')$ where $G'$ is the graph induced by $V\setminus I_1$ must be ve-dominated by $I_2$ and both endpoints of those edges have a neighbour in $I_2$. Since $I_1$ and $I_2$ are disjoint, every edge in $V\setminus I_1$ is $2$-ve dominated by $I_1\cup I_2$. Now consider an edge in the graph induced by $V\setminus(I_1\cup I_2)$. Clearly, that edge is $3$-ve dominated by $I_1\cup I_2\cup I_3$ using the similar arguments. Now, consider an edge incident to some vertex in $I_1$. Since these edges have an endpoint in $V\setminus I_2$, those endpoints must have a neighbour in $I_2$. Therefore, every edge incident to $I_1$ is also $2$-ve dominated. Hence we have $|N_G[e]\cap L|\geq 2$, for every $e\in E$.
	
	Now, suppose that there is a pair of edge $e,e'$ in $E$ such that $|(N_G[e]\cup N_G[e'])\cap L|<3$. Clearly, $e$ and $e'$ has no endpoint in $V\setminus(I_1\cup I_2)$. So, let us assume that $e$ and $e'$ have one endpoint in $I_1$ and one endpoint in $I_2$. The only way $|(N_G[e]\cup N_G[e'])\cap L|<3$ happens if $e$ and $e'$ are incident to same pair of vertices in $I_1\cup I_2$, which is not possible. Therefore, no such pair of edges exists in $G$. Hence, for every pair of edges $e,e'$ in $E$, we have $|(N_G[e]\cup N_G[e'])\cap L|\geq 3$. Thus, $L$ is a liar's ve-dominating set of $G$.
\end{proof}

Now, we give an algorithm to construct the desired $2$-separated collection. At the beginning, we set $r=1$. We keep on increasing the value of $r$ until $|L(N_G^{r+6}[v])|>\rho |L(N_G^r[v])|$ holds. Let $r'$ be the smallest $r$ violating the condition. We set $\mathcal{S}_i=N_G^{r'}[v]$ and $\mathcal{T}_i=N_G^{r'+6}[v]$. After that, we update $V_i$ by $V_i\setminus N_G^{r'+2}[v]$. Also, we remove all the singletone vertices in $V_i\setminus N_G^{r'+2}[v]$. We repeat the steps until $V_i$ is non-empty. Suppose that there are $k$ number of iterations of the algorithm. Note that, the collection $\mathcal{S}=\{\mathcal{S}_1,\mathcal{S}_2,\ldots,\mathcal{S}_k\}$ is a $2$-separated collection. Thus, the algorithm returns a $2$-separated collection $\mathcal{S}=\{\mathcal{S}_1,\mathcal{S}_2,\ldots,\mathcal{S}_k\}$ and it's corresponding collection $\mathcal{T}=\{\mathcal{T}_1,\mathcal{T}_2,\ldots,\mathcal{T}_k\}$. The pseudocode is given in Algorithm \ref{Algo:PTAS}.

\begin{algorithm}
	\textbf{Input:} A unit disk graph $G=(V,E)$ and an arbitrary small $\epsilon>0$.\\
	\textbf{Output:} The sets $\mathcal{S}=\{\mathcal{S}_1,\mathcal{S}_2,\ldots,\mathcal{S}_k\}$ and $\mathcal{T}=\{\mathcal{T}_1,\mathcal{T}_2,\ldots,\mathcal{T}_k\}$.
	\caption{A $2$-separated collection}
	\label{Algo:PTAS}
	\begin{algorithmic}[1]
		\State  $V_0=V$;
		\State $i=0$ and $\rho=1+\epsilon$;
		\While{$V_i\neq \emptyset$}
		\State Choose a vertex $v\in V_i$;
		\State $N^0_G[v]=\{v\}$ and $r=1$;
		\While{$|L(N_G^{r+6}[v])|>\rho |L(N^r[v])|$} \Comment{Call Algorithm \ref{Algo:liarindepent}}
		\State $r=r+1$;
		\EndWhile
		\State $r'=r$;
		\State $i=i+1$;
		\State $\mathcal{S}_i=N_G^{r'}[v]$ and $\mathcal{T}_i=N_G^{r'+6}[v]$;
		\State $V_i=V_{i-1}\setminus N_G^{r'+2}[v]$;
		\For{every isolated vertex $v'$ in $V_i$}
		\State $V_i=V_i\setminus \{v'\}$;
		\EndFor
		\EndWhile
		\State \Return $\mathcal{S}=\{\mathcal{S}_1,\mathcal{S}_2,\ldots,\mathcal{S}_k\}$ and $\mathcal{T}=\{\mathcal{T}_1,\mathcal{T}_2,\ldots,\mathcal{T}_k\}$;
	\end{algorithmic}
\end{algorithm}

First, we prove that the $L(N_G^{r}[v])$ is a liar's ve-dominating set of $N^{r}_G[v]$ and the size is bounded by $O(r^2)$.
\begin{lem}\label{lem:liardomunitdiscbounded}
	$L(N_G^r[v])$ is a liar's ve-dominating set of $N_G^r[v]$ in a unit disk graph $G$ and $|L(N^r_G[v])|\leq O(r^2)$.
\end{lem}
\begin{proof}
	Since, we compute $L(N_G^r[v])$ using the Algorithm \ref{Algo:liarindepent}, from Lemma \ref{lem:liarindependent} it follows that $L(N_G^r[v])$ is a liar's ve-dominating set of $N_G^r[v]$. Also, $|L(N_G^r[v])|\leq |I_1|+|I_2|+|I_3|$ where $I_1$, $I_2$ and $I_3$ are independent sets obtained by Algorithm \ref{Algo:liarindepent}. Since $G$ is a unit disk graph the size of the maximal independent set in $N_G^r[v]$ is bounded by the number of non-intersecting unit disks packed in a larger disk of radius $r+\frac{1}{2}$ centered at $v$. Therefore, $|I_i|\leq \frac{\pi(2r+1)^2}{\pi (1)^2}$ for every $i=1,2,3$. Thus, we have $|L(N_G^r[v])|\leq 3\frac{\pi(2r+1)^2}{\pi (1)^2}=3(2r+1)^2$ and hence the bound. 
\end{proof}

In the next lemma, we show that there is an $r$ violating the condition $|L(N_G^{r+6}[v])|>\rho |L(N^r[v])|$.

\begin{lem}\label{lem:smallestrviolatingcondn}
	There exists an $r$ violating the condition $|L(N^{r+6}_G[v])|>\rho |L(N^r_G[v])|$, for $\rho =1+\epsilon$.
\end{lem}

\begin{proof}
	We prove this using contradiction. Let us assume that there exists a vertex $v$ such that $|L(N^{r+6}_G[v])|>\rho L(N^r_G[v])$ for every $r=1,2,\ldots$. From Lemma \ref{lem:liardomunitdiscbounded}, we have $|L(N^r_G[v])|\leq 3(2r+1)^2$. Observe that $|L(N^1_G[v])|\geq 2$ and $|L(N^2_G[v])|\geq 3$. Now, if $r=6k$, we have $3(2r+1)^2\geq |L(N^r_G[v])|>\rho |L(N^{r-6}_G[v])|>\ldots>\rho^{\frac{r}{6}}|L(N^2_G[v])|\geq 3\rho^{\frac{r}{6}}$. Also, if $r=6k+s$ where $1\leq s\leq 5$, we have $3(2r+1)^2\geq |L(N^r_G[v])|>\rho |L(N^{r-6}_G[v])|>\ldots>\rho^{\frac{r-5}{6}}|L(N^1_G[v])|\geq 2\rho^{\frac{r-5}{6}}$. Hence, we have
	\begin{align}
		12(2r+1)^{12}>& \rho^{r},~if~ r=6k\\
		>&\rho^{r-5},~ if~  r=6k+s, ~where~ 1\leq s\leq 5
	\end{align}
	The left hand side of the above inequalities is a polynomial function in $r$ but the right hand side is an exponential function in $r$. For arbitrary large $r$ these inequalities do not hold which is a contradiction. Hence, there is a value of $r$, say $r'$ for which the condition is violated.  
\end{proof}

Next we show that the smallest $r$ violating the condition is bounded by $O(\frac{1}{\epsilon}\log\frac{1}{\epsilon})$.

\begin{lem}\label{lem:boundonrptas}
	The smallest $r$ violating the condition is bounded by $O(\frac{1}{\epsilon}\log\frac{1}{\epsilon})$.
\end{lem}

\begin{proof}
	Let $r'$ be the smallest $r$ violating the condition. We show that $r'=O(\frac{1}{\epsilon}\log\frac{1}{\epsilon})$ using the following inequalities:
	
	\noindent$(i)$ $\log(1+\epsilon)> \frac{\epsilon}{2}$ for $0<\epsilon<1$
	
	\noindent$(ii)$ $\log x<x$ for $x>1$
	
	\noindent$(iii)$ $\log(\frac{1}{\epsilon})\geq 1$ for $\epsilon\leq \frac{1}{10}$
	
	Suppose that $x=\frac{c}{\epsilon}\log(\frac{1}{\epsilon})$. Also, consider the inequality $12(2x+1)^{12}\leq 12(3x)^{12}\leq (1+\epsilon)^x$. Taking logarithm on both sides of this inequality we have $\log 12+12\log 3x\leq x\log(1+\epsilon)$. Using inequality $(1)$, we have $\frac{\log 12+12\log 3x}{x}\leq \frac{\epsilon}{2}$. Therefore, to show the bound it is sufficient to show that the above-mentioned inequality holds. Using inequality $(ii)$ and $(iii)$, we can show that the choice of $x$ satisfy the inequality if $\log 3c+1<\frac{c}{24}$, for any constant $c$.  
\end{proof}

From Lemma \ref{lem:liardomunitdiscbounded} and Lemma \ref{lem:boundonrptas}, it implies that the size of a liar's ve dominating set of $\mathcal{T}_i$ is bounded by a constant. In the following lemma, we show that $L^*(\mathcal{T}_i)$ can be computed in polynomial time.

\begin{lem}\label{lem:PTASOptimalinpolynomialtime}
	For a given $v\in V$, $L^*(\mathcal{T}_i)$ can be computed in polynomial time.
\end{lem}

\begin{proof}
	For a given $v\in V$, we have $\mathcal{T}_i\subseteq N_G^{r+3}[v]$. Consider the graph induced by $\mathcal{T}_i$. By Lemma \ref{lem:liardomunitdiscbounded}, the size of $L^*(\mathcal{T}_i)$ is bounded by $O(r^2)$. Therefore, $L^*(\mathcal{T}_i)$ can be computed by considering every possible subset of $V$ of size at most $r^2$ in $G$ and checking whether it is a liar's ve-dominating set of $\mathcal{T}_i$ or not. There are $O{n\choose r^2}$ subsets of $V$ of size at most $r^2$ and checking a subset of $V$ is liar's ve-dominating set of $\mathcal{T}_i$ or not can be done in polynomial time. In Lemma \ref{lem:boundonrptas}, we have shown that $r$ is bounded by a constant. Therefore, the $L^*(\mathcal{T}_i)$ can be determined in polynomial time.
\end{proof}

Next, we show that $\bigcup\limits_{i=1}^{k}L(\mathcal{T}_i)$ is a liar's ve-dominating set of $G$.

\begin{lem}
	For the collection of neighbourhoods $\{\mathcal{T}_1,\mathcal{T}_2,\ldots,\mathcal{T}_k\}$ given by Algorithm \ref{Algo:PTAS}, $L=\bigcup\limits_{i=1}^{k}L(\mathcal{T}_i)$ is a liar's ve-dominating set in $G$.
\end{lem}

\begin{proof}
	First, we show that for every edge $e\in E$, $|N_G[e]\cap L|\geq 2$. Observe that, every edge belongs to some $\mathcal{T}_i$ by the construction given in Algorithm \ref{Algo:PTAS}. By Lemma \ref{lem:liardomunitdiscbounded}, $L(\mathcal{T}_i)$ is a liar's ve-dominating set of $\mathcal{T}_i$. This implies that $|N_G[e]\cap L(\mathcal{T}_i)|\geq 2$ for some edge $e$ in the graph induced by $\mathcal{T}_i$. Therefore $|N_G[e]\cap L|\geq 2$, for any edge $e\in E$.
	
	Now, we show that for any pair of edges $e,e'$, $|(N_G[e]\cup N_G[e'])\cap L|\geq 3$. If $e,e'$ is in the graph induced by same $\mathcal{T}_i$, then $|(N_G[e]\cup N_G[e'])\cap L|\geq 3$ is satisfied because $L(\mathcal{T}_i)$ is a liar's ve-dominating set of $\mathcal{T}_i$. If $e$ is an edge of the graph induced by $\mathcal{T}_i$ and $e'$ is an edge of the graph induced by $\mathcal{T}_j$ with $i\neq j$ and $N_G[e]\cap L(\mathcal{T}_i)\neq  N_G[e']\cap \mathcal{T}_j$, then obviously we have $|(N_G[e]\cup N_G[e'])\cap L|\geq 3$. Now suppose that, $e$ is in some graph induced by $\mathcal{T}_i$ and $e'$ is in some graph induced by $\mathcal{T}_j$ with $i\neq j$ and $N_G[e]\cap L(\mathcal{T}_i)=  N_G[e']\cap L(\mathcal{T}_j)$. Observe that if $|N_G[e]\cap L(\mathcal{T}_i)|\geq 3$, then the pair $e,e'$ satisfy $|(N_G[e]\cup N_G[e'])\cap L|\geq 3$. Therefore, let $|N_G[e]\cap L|=|N_G[e']\cap L|=2$. By the construction of the collection $\{\mathcal{T}_1,\mathcal{T}_2,\ldots,\mathcal{T}_k\}$, the edges which belongs to different subgraphs $\mathcal{T}_i$ and $\mathcal{T}_j$ but share some neighbours must be in the $\mathcal{T}_i\cap\mathcal{T}_j$. Therefore, for those edges we have $|(N_G[e]\cup N_G[e'])\cap L|\geq 3$. Hence, $L$ is a liar's ve-dominating set in $G$.  
\end{proof}

\begin{coro}\label{cor:optimal=unionofsuboptimals}
	For the collection of neighbourhoods $\{\mathcal{T}_1,\mathcal{T}_2,\ldots,\mathcal{T}_k\}$, $L^*=\bigcup\limits_{i=1}^{k}L^*(\mathcal{T}_i)$ is a liar's ve-dominating set in $G$.
\end{coro}

\begin{theo}
	Given a unit disc graph $G=(V,E)$ and any $\epsilon>0$, the liar's ve-dominating set can be approximated within a factor of $(1+\epsilon)$, where the corresponding approximation algorithm has running time $n^{O(c^2)}$ and $c=O(\frac{1}{\epsilon}\log \frac{1}{\epsilon})$.
\end{theo}

\begin{proof}
	The proof follows from Lemmas \ref{lem:unitdiscsecond}, \ref{lem:smallestrviolatingcondn}, \ref{lem:boundonrptas}, \ref{lem:PTASOptimalinpolynomialtime} and Corollary \ref{cor:optimal=unionofsuboptimals}.
\end{proof}

\section{Conclusion}
In this article, we have proved that the \textsc{DecideLVEDP} is NP-complete in unit disk graphs. Also, we designed a PTAS for \textsc{MinLVEDP} in the unit disk graph. It is interesting to study the liar's ve-domination problem for other graph classes of bipartite and chordal graphs. 
\bibliographystyle{plain}
\bibliography{VEDom_bib}

\end{document}